\documentclass[12pt]{article}
\usepackage{graphics}

\date{}
\def\qed{{\hfill\rule{1.2ex}{1.2ex}}}

\newenvironment{Definition}{\begingroup\noindent{\bf Definition:}}{\endgroup}
\newtheorem{Lemma}{Lemma}

\def\:{\mathop{:}\nolimits}

\def\R{{\bf R}}
\def\Q{{\bf Q}}
\def\N{{\bf N}}
\def\Z{{\bf Z}}

\def\Id{\mathop{\hbox{\rm Id}}\nolimits}

\def\Sup{\mathop{\rm Sup\,}\nolimits}

\author{Norbert A'Campo}
\title{A natural construction for the real numbers}

\begin{document}

\maketitle
\noindent
We propose  a new construction of the real number system, that
is built directly upon the additive group of integers 
and has its roots in the definition due to 
Henri Poincar\'e [P, pages 230--233,]
of the rotation number of an orientation preserving
homeomorphism of the circle. The definitions of 
addition, multiplication and comparison 
of real numbers are very natural in our setting. 
The proposed definition 
of real numbers is illustrated with examples that are irrational, roots of an integral polynomial 
equation, but not expressible by radicals, or not root of
an integral polynomial equation.
I thank Sebastian Baader, Etienne Ghys and Domingo Toledo 
for stimulating
discussions.

\noindent
{\bf Slopes and definition of the real numbers.} 

\noindent
Let $(\Z,+)$ be 
the set of integers together
with the arithmetic operation
of addition. 
The basic objects in our construction are slopes. 
A {\it slope} is by definition a map
$\lambda:\Z \to \Z$, with the property that the set
$\{\lambda(m+n)-\lambda(m)-\lambda(n) \mid m,n \in \Z\}$ 
is finite. Two slopes $\lambda, \lambda'$
are {\it equivalent} if the set $\{\lambda(n)-\lambda'(n) \mid n\in \Z\}$
is finite.

\begin{Definition}
A {\it real number} is an equivalence class of slopes.
\end{Definition}

Let $\R$ denote the set of
real numbers.
For $j \in \Z$, let $\bar{j}:\Z \to \Z$ be the map $\bar{j}(n):=nj$.
The linear map $\bar{j}$ is a slope for which the expression
$\bar{j}(n+m)-\bar{j}(n)-\bar{j}(m)$ takes only the value $0$. We
identify an integer $j \in \Z$ with the real number
represented by the slope $\bar{j}$. After this identification the set of
integers $\Z$ becomes a subset of the set of real numbers $\R$.
We see
that among the real numbers the integers appear as
those real numbers, which are representable by a
linear slope.

For $p, q \in \Z, q > 0, $ let the map $\phi:\Z \to \Z$ be defined by
$\phi(n):= \min \{k \in \N \mid  qk \geq pn\}, n \in \N, n>0$ and by 
$\phi(-n)=-\phi(n)$ for $n \in \Z, n \leq 0$.
The map $\phi$ is a slope representing  the rational number
$p/q$, i.e.
the slope $\phi$ represents a real number which is a solution of
the equation $qx=p$. This will become clear, when we have defined multiplication of real numbers in our setting. 
As for the integers we identify the set of
rational numbers $\Q$ with a subset of $\R$. One can characterize the rational
numbers as those real numbers, which are representable by a slope $\lambda$,
such that for some integer $q >0$ the map
$n \in \Z \mapsto \lambda(qn) \in \Z$ is linear.

We now
define the basic arithmetic operations such as {\it addition} and
{\it multiplication}
of real numbers.

Let $a, b \in \R$ be real numbers. Let $\alpha, \beta $ be
slopes representing the real numbers $a$ and $b$. The map
$\alpha+\beta: \Z \to \Z$, which is defined by
$(\alpha+\beta)(n):=\alpha(n)+\beta(n)$, is again a 
slope and its equivalence class is 
independent of the choice of representatives
$\alpha, \beta $ for $a,b$. We define the 
{\it sum} $a+b \in \R$ of $a, b \in \R$ 
as the equivalence class of the slope $\alpha+\beta: \Z \to \Z$.

The composition $\alpha \circ \beta:\Z \to \Z$ is
again a slope, and we
define the {\it product} $ab \in \R$ as the equivalence class of the composition
$\alpha \circ \beta:\Z \to \Z$.

The consistency of this definition follows from the following lemma.

\begin{Lemma} Let the slopes 
$\alpha, \alpha'$ represent $a\in \R$ and the slopes $ \beta, \beta'$ represent $b\in \R$.
Then the compositions $\alpha \circ \beta$ and $\alpha' \circ \beta'$ are equivalent
slopes.
\end{Lemma}

\noindent
{\bf Proof.\enspace}
We first show that the map $\alpha \circ \beta$ is a slope. Let $E_{\alpha}$
and $E_{\beta}$ be finite subsets in $\Z$, such that
$\alpha(n+m)-\alpha(n)-\alpha(m) \in E_{\alpha}$ and 
$\beta(n+m)-\beta(n)-\beta(m) \in E_{\beta}$ for $n,m \in \Z$.
Hence, for $n,m \in \Z$ there exist $u,u' \in E_{\alpha},v\in E_{\beta}$ with
$$
\alpha \circ \beta(n)+\alpha \circ \beta(m)-\alpha \circ \beta(n+m)=
$$
$$
\alpha(\beta(n)+\beta(m))+u-\alpha(\beta(n)+\beta(m)-v)=
$$
$$
\alpha(\beta(n)+\beta(m))+u-(\alpha(\beta(n)+\beta(m))+\alpha(-v)-u')=u-\alpha(-v)-u'.
$$
We conclude that the expression 
$\alpha \circ \beta(n)+\alpha \circ \beta(m) -\alpha \circ \beta(n+m), 
n,m \in \Z,$ 
takes  its values in a finite set. Hence, the map $\alpha \circ \beta$
and, with the same justification, also the map  
$\alpha' \circ \beta'$ are slopes.

Let $E_{\alpha,\alpha'}$ and $E_{\beta,\beta'}$ be finite sets such 
that we have
$\alpha(n)-\alpha'(n) \in E_{\alpha,\alpha'}$ and 
$\beta(n)-\beta'(n) \in E_{\beta,\beta'}$ for $n \in \Z$.
Hence, for $n \in \Z$ there exist
$r \in E_{\alpha,\alpha'}, s \in E_{\beta,\beta'}$ and 
$u \in E_{\alpha}$ with
$$
\alpha \circ \beta(n)-\alpha' \circ \beta'(n)=
\alpha(\beta'(n)-s)-(\alpha(\beta'(n))+r)=
$$
$$
\alpha(\beta'(n))+\alpha(-s)-u-(\alpha(\beta'(n))+r)=
\alpha(-s)-r-u.
$$
We conclude that the expression 
$\alpha \circ \beta(n)-\alpha' \circ \beta'(n), n \in \Z,$ 
takes its values in a finite set.
Hence, the slopes $\alpha \circ \beta$ and $\alpha' \circ \beta'$
are equivalent.
\qed

 Let
 $\N:=\{n \in \Z \mid n \geq 0 \}$ be the set of natural numbers.
We call a slope $\lambda$ {\it positive}, if the set
$\{\lambda(n) \mid  n\in \N, \lambda(n) \leq 0 \}$ is finite, while the
set $\{\lambda(n), n\in \Z\}$, is infinite. A real number $a$ is 
{\it positive}, if its representing slopes are positive.

From this definition of positivity we obtain the ordering of the
real numbers as usual in the following way.
If $a$ is positive, we  say that $a >  0$ and $0 < a$ hold. The
real number $a$ is defined to be {\it less} then the
real number $b$ if there exists a
positive real number $t$ with $b=a+t$. If  $a$ is
less then $b$, we say that  $a < b$ holds.

We illustrate  the definitions by examples before
stating and verifying that the set $\R$ with the
addition $+$ , multiplication $\cdot$ and order relation $<$ satisfies
all the axioms of the real numbers, i.e. of a complete totally ordered
archimedean field.

A map $f:\Z \to \Z$ is called {\it odd} if for all $n \in \Z$ the property
$f(-n)=-f(n)$ holds. An odd map $f:\Z \to \Z$ is determined by its
restriction to $\N_+:=\{n \in \Z \mid n > 0 \}$. Let $\lambda$ be an arbitrary
slope. Then the  map $\kappa:\Z \to \Z$ with $\kappa(0)=0$ defined by 
$\kappa(n):=\lambda(n), n > 0,$ and by $\kappa(n):=-\lambda(-n), n < 0,$
is an odd slope, which is equivalent to the slope $\lambda$. So every
real number can be represented by an odd slope.
In
order to verify that an odd map $\gamma:\Z \to \Z$ is a slope, it suffices
to check that the set $\{\gamma(n+m)-\gamma(n)-\gamma(m) \mid n,m \in \N_+\}$
is finite.

We will construct a slope that represents the number $\sqrt{2}$. 
Let $\rho:\Z \to \Z$ be the odd map defined by 
$\rho(n)=\min \{ k \in \N \mid 2n^2 \leq k^2 \}, n \in \N_+$. We have
for $n \in \N_+$  the inequalities:
$n \leq \rho(n) \leq 2n, 2n^2 \leq \rho(n)^2, (\rho(n)-1)^2 \leq 2n^2$. Hence
$2n^2 \leq \rho(n)^2 \leq 2n^2+2\rho(n)-1 \leq 2(n+1)^2$. For $n,m \in \N_+$ we deduce
$2nm \leq \rho(n) \rho(m) \leq 2(n+1)(m+1)$.
The map $\rho$ is a slope, since for $n,m \in \N_+$ we estimate
$$
x:=(-\rho(n+m)+\rho(n)+\rho(m))(\rho(n+m)+\rho(n)+\rho(m))= 
$$
$$
-\rho(n+m)^2+\rho(n)^2+\rho(m)^2+2\rho(n)\rho(m)
$$
by
$$
-4n-4m-2=-2(n+m+1)^2+2n^2+2m^2+4nm \leq x  \leq 
$$
$$
-2(n+m)^2+2(n+1)^2+2(m+1)^2+4(n+1)(m+1)=8m+8n+8
$$
and with $\rho(n+m)+\rho(n)+\rho(m) \geq n+m+1$, we conclude
$|\rho(n+m)-\rho(n)-\rho(m)| \leq 8$.
The equivalence class of $\rho$ is a
positive real number $a$ satisfying $a^2=2$. Indeed, the number 
$a^2$ is represented by the composition 
$\rho \circ \rho$. We have for $n \in \N_+$ the inequalities
$4n^2\leq 2\rho(n)^2\leq \rho(\rho(n))^2\leq 2(\rho(n)+1)^2\leq4n^2+8n+2\leq 4(n+1)^2$
showing $2n \leq \rho(\rho(n)) \leq 2n+2$, which means that 
the slopes $\rho \circ \rho$ and $\bar{2}$ are equivalent and 
represent the integer $2$. Hence $\rho$
represents the square root $\sqrt{2}$ of $2$, which is the length of the
diagonal of a unit square and 
can not be represented as a fraction ${p \over q}$, see [F,V].

We will now construct a real number that is a 
root of the polynomial $p(x):=x^5+x-3$. Let the odd map $\alpha:\Z \to \Z$ 
be defined by 
$\alpha(n):=\min \{ k \in \N \mid 3n^5 \leq k^5 +n^4k \}=\min\{k \in \N \mid p({k \over n}) \geq 0\}, n \in \N_+.$ 
Observe that
that 
$p({{\alpha(n)-1} \over n})\leq 0 \leq p({{\alpha(n)} \over n})$ and ${{\alpha(n)} \over n}\leq {6\over 5}$ hold and that
$p({{\alpha(n)} \over n}) =p({{\alpha(n)-1} \over n}+{1 \over n})\leq {6131 \over {125n}} \leq {50 \over n}$ follows.
The map $\alpha$ is a slope and represents the real 
root $a$ of the equation
$x^5+x-3=0$, which can not be represented by a compound radical expression
after the work of Paolo Ruffini (1762-1822) and of Niels Henrik Abel 
(1802--1829), see [A,R,S]. 
We show that $\alpha$ is a slope. For $n,m \in \N_+$ define the 
rational numbers
$a_-:=\max\{{{\alpha(n)-1} \over n}, {{\alpha(m)-1} \over m},
{{\alpha(m+n)-1} \over {m+n}}\}$ and $a_+:= \min \{{{\alpha(n)} \over n}, {{\alpha(m)} \over m},
{{\alpha(m+n)} \over {m+n}}\}$. From the 
monotonicity of $p$ and the definition
of $\alpha$ we deduce $p(a_-) \leq 0, p(a_+) \geq 0$ 
and $a_- \leq a_+$. Let $A$
be any rational number with $a_- \leq A \leq a_+$. We have
the inequalities $|\alpha(n)-nA| \leq 1,|\alpha(m)-mA| \leq 1$ and 
$|\alpha(m+n)-(m+n)A| \leq 1$, that show
$|\alpha(m+n)-\alpha(m)-\alpha(n)| \leq 3$.
Let $a$ be the real number that is represented by the slope $\alpha$.
We show $p(a)=0$ by showing that the slope 
$\alpha^{\circ 5}+\alpha-\bar{3}$ is bounded. Here we have used 
the notation
$\alpha^{\circ e}$ for the $e$-th iterate, $e\in \N$, of $\alpha$.
It is not at all easy to handle directly 
iterates of slopes. The following
estimate helps out and is proved by induction upon the exponent $e \in \N_+$
$$
|n^{e-1}\alpha^{\circ e}(n)-\alpha(n)^e| \leq n^{e-1} (1+|\alpha(1)|+S_{\alpha})^{e-1}, n \in \N_+,
$$
where the quantity  
$S_{\alpha}:=\max\{|\alpha(u+v)-\alpha(u)-\alpha(v)|, u,v \in \Z\}$ 
measures the non-linearity of the slope $\alpha$. 
Note also $|\alpha(n)|\leq |n|(|\alpha(1)|+S_{\alpha})$.
It follows that the slope $\alpha^{\circ e}$ is equivalent to the odd slope
defined by $n\in \N_+ \mapsto [{{\alpha(n)^e} \over {n^{e-1}}}]$. 
So the slope
$\alpha^{\circ 5}+\alpha-\bar{3}$ is equivalent to the odd slope 
$\epsilon$ defined by $n\in \N_+ \mapsto [{{\alpha(n)^5} \over {n^{4}}}]+\alpha(n)-3n$. 
The slope $\epsilon$ is bounded, since 
for $n \in \N_+$ we have  $o \leq \epsilon(n)=np({{\alpha(n)} \over n})\leq 50$. 
Here $[?]$ is the Gaussian integral part bracket 
$[x]:=\max\{k \in \Z \mid k \leq x \}, x \in \R$.

Let $\beta:\Z \to \Z$ be the odd map with 
$\beta(n):=\# \{ (p,q) \in \Z \times \Z \mid p^2+q^2 \leq n \}, n > 0$.
Unit squares in the plane 
with centers at the lattice points 
$(p,q)\in \Z \times \Z, p^2+q^2 \leq n,$
cover the disk of radius $\sqrt{n}-{1 \over 2}\sqrt{2}$ and are contained
in the disk of radius $\sqrt{n}+{1 \over 2}\sqrt{2}$.
Hence $|\beta(n)-n\pi| \leq 2\sqrt{2}\sqrt{n}, n \in \N$.
It follows that the odd map $\bar{\beta}$ defined by 
$n\in N_+ \mapsto [{{\beta(n^2)} \over n}]$ is a slope.  The 
slope $\bar{\beta}$ represents the area $\pi$ 
of the unit disk in the plane, see
Chap. $1$ in {\it Anschauliche Geometrie} of 
David Hilbert \& S. Cohn-Vossen [H-V]. 
Johann Heinrich Lambert (1728-1777) has proved in his communication 
of 1761 to the Academy in Berlin  
that the number $\pi$ is not root of a polynomial equation with
integral coefficients, see [Le].
The map $\beta$ is not a slope since the quantity $s(n):=\beta(n)-\beta(n-1)-\beta(1)$ is not 
bounded. One has for instance $s(5^u)=4u-1, u\in \N_+$.

The number 
$e$ appeared in the
sixteenth century, when it was noticed that the expression $(1+{1 \over n})^n$
for compound interest increases with $n$ to a certain value 
$2.7182818\cdots$, see the book ``{\bf e} {\it The
Story of a Number}'' 
by Eli Maor [M]. The number $e$ became of central importance
in Mathematics since its interpretations in Geometry and Analysis
by Gr\'egoire de Saint-Vincent (1584-1667). It is not obvious to define
the number $e$ with a slope. We use the solution to 
a problem, see [D],  of Jakob Steiner 
(1796-1863) and define for $n \in \N, n > 0,$ the integer
$\epsilon(n)$ to be the natural number $k, k>0,$ such that the expression 
$({k \over n})^{{n \over k}}$ takes its maximal value.
The corresponding 
odd function $\epsilon: \Z \to \Z$ is a slope representing the number $e$.

The classical construction of the system of real numbers
is based on
Dedekind cuts or  on Cauchy sequences $(r_n)_{n \in N}$
of rational numbers. The present construction by slopes is related
to the classical ones as follows: To a slope
$\lambda$ corresponds a Dedekind cut $(A,B)$
by setting
$A:=\{{p \over q} \in \Q \mid \bar{p}  \leq \lambda \circ \bar{q} \}$
and $B:=\{{p \over q} \in \Q \mid \lambda \circ \bar{q}  \leq \bar{p} \}$ and
also a Cauchy sequence $(r_n)_{n \in N}$ by setting $r_n:={\lambda(n+1) \over {n+1}}$.

\noindent
{\bf Well adjusted slopes.\enspace}

We call a slope $\lambda$ {\it well adjusted\/}
if it is odd and satisfies the inequalities
$-1 \leq \lambda(m+n)-\lambda(m)-\lambda(n) \leq 1, n,m \in \Z$.
One can say that a well adjusted slope need not be a linear map 
from $\Z$
to $\Z$, but deviates as little as possible from being linear.  
Each slope is equivalent to a 
well adjusted slope, as  shows the concentration Lemma below. So
in particular,
a real number can be represented by a well adjusted slope.

\noindent 
For integers $p, q, q\not=0,$ the result of  
optimal euclidean division
of $p$ by $q$ will be denoted by $p:q$. The {\it optimal euclidean division} 
is the  integer $r:=p\: q \in \Z$ that satisfies the inequalities
$2p-|q| \leq 2qr < 2p+|q|$, where 
$|q|:=\max\{q,-q\}$ is the absolute value of $q$. For instance $4\: 7=1$ but
$3\: 7=0$. If $p/q, p,q \in \Z, q\not=0,$ denotes the fraction, then we have
$|p/q-p\: q| \leq 1/2$. For the optimal euclidean division, we have 

\begin{Lemma} Let $q \in \N_+$ and $a,b,c \in \Z$ be such that $-q \leq a-b-c \leq q$. Then we have
$-1 \leq a\: 3q-b\: 3q-c\: 3q \leq 1$.
\end{Lemma}

\noindent
{\bf Proof.\enspace}
The integer $a\: 3q-b\: 3q-c\: 3q$ differs from $0$ by at most
$1/2+1/2+1/2+|a/3q-b/3q-c/3q|\leq 3/2+1/3=11/6$. 
Hence we have $-1 \leq a\: 3q-b\: 3q-c\: 3q \leq 1$, since $11/6<2$.
\qed

\begin{Lemma} Let $n,m \in \N_+$ and $c \in \Z$. Then we have
$$
-1 \leq c\: m(n+m)-c\: n(n+m)-c\: nm \leq 1.
$$
\end{Lemma}

\noindent
{\bf Proof.\enspace}
The integer $c\: m(n+m)-c\: n(n+m)-c\: nm$ differs from 
$c/m(n+m)-c/n(n+m)-c/nm=0$ 
by at most
$1/2+1/2+1/2=3/2$, hence $-1 \leq c\: m(n+m)-c\: n(n+m)-c\: nm \leq 1$,
since $3/2<2$.
\qed

\begin{Lemma} [Concentration Lemma] 
Let $\lambda$ be a slope. Let $s \in \N_+$ be such that
for all $n,m \in \Z$ we have
$-s \leq \lambda(m+n)-\lambda(m)-\lambda(n) \leq s$.
Let $\lambda':\Z \to \Z$ be defined by
$\lambda'(n):=\lambda(3sn)\: 3s, n \in \Z$.
Then the map $\lambda'$ is a well adjusted slope, which is  
equivalent to the slope
$\lambda$.
\end{Lemma}

\noindent
{\bf Proof.\enspace}
By induction on $t \in \N_+,$  we prove
$-s(t-1) \leq \lambda(tn)-t\lambda(n) \leq s(t-1)$.
For $t=3s$ we get 
$-s(3s-1) \leq \lambda(3sn)-3s\lambda(n) \leq s(3s-1)$ and hence
$-s \leq \lambda'(n)-\lambda(n) \leq s$, which shows the equivalence of 
$\lambda$ and $\lambda'$. From
$-s \leq \lambda(3sn+3sm)-\lambda(3sn)-\lambda(3sm) \leq s$  we
deduce $-1 \leq \lambda'(n+m)-\lambda'(n)-\lambda'(m) \leq 1$.
\qed

A well adjusted slope $\lambda$ has the following properties:\newline
\noindent -  $|\lambda(n+1)-\lambda(n)| \leq |\lambda(1)|+1,$\newline
\noindent - if for some $k \in \N_+$ we have $\lambda(k) > 1$ 
(or $\lambda(k) < -1$), then we have
for any $n \in \N_+$ the inequality $\lambda(n) \geq -1+n\: k$ 
(or $\lambda(n) \leq +1-n\: k$),\newline
\noindent - if for some $k \in \Z$ we have $\lambda(k) > 1$, then for $v \in \Z$
the set $\{n \in \Z \mid \lambda(n)=v \}$ is finite and has fewer then $k+1$
elements,\newline
\noindent - if for some $k \in \Z$ we have $\lambda(k) > 1$, then for any $v \in \Z$, there exists $n \in \Z$ with 
$|v-\lambda(n)| \leq |\lambda(1)|+1,$\newline
\noindent - the real number $x$ represented by $\lambda$ satisfies
$x > 0$ if and only if there exists $a \in \N$ with $\lambda(a)>1$,\newline
\noindent - let $y$ be a real number represented by a well adjusted slope $\kappa$.
We have $x > y$ if and only if there exists $n \in \N_+$ with 
$\lambda(n)>2+\kappa(n)$.

From the above lemma and remarks we obtain:

\begin{Lemma} Let $\lambda$ be a slope. If 
$\lambda$ takes infinitely many values,
then there exist $b,B \in \N_+$ such that the following inequalities hold:
$$
|\lambda(n+k)-\lambda(n)| \leq kb, n \in \Z, k \in \N,
$$
$$
|\lambda(n+kB)-\lambda(n)| \geq k, n \in \Z, k \in \N.
$$
In particular, the slope $\lambda$ takes each value at most $2B-1$ times.
\qed
\end{Lemma}

\noindent
{\bf The axioms.\enspace}

\noindent
We now state, partially in abbreviated form,  the axioms for a complete totally ordered field, that are 
satisfied by the quadruple
$(\R,+, \cdot ,<)$. The presentation of the axioms is slightly redundant.

\begin{enumerate}
\item The pair $(\R,+)$ is an abelian group.

\item The triple $(\R,+, \cdot)$ is a field.

\item The quadruple $(\R,+,\cdot,<)$ is an archimedean, complete, totally 
ordered field. 

 Complete ordered field, i.e.
 \begin{itemize}
 \item for any non-empty subset $T$ bounded from above in $ \R$ there exists 
 a least upper bound in $ \R$ called the {\it supremum} of $T$. It 
 will be denoted by $\Sup T$.
 \end{itemize}
  
Archimedian ordered field, i.e.
  \begin{itemize}
  \item for $a \in \R,a > 0$ and $A \in \R$ there exists a $N \in \N$ 
  such that $Na > A$.
  \end{itemize}
  \end{enumerate}

We now begin the verification of the axioms for the system of real
numbers, that we have introduced above. We leave out those verifications, that are straightforward and can be done without
using well adjusted slopes as representatives.

The addition $+$ of integers makes $Z$ into an abelian group $(Z,+)$. It
follows easily that $(\R,+)$ is also an abelian group. 

The triple $(\R,+,\cdot)$ is a skew field, i.e. verifies all the field
axioms. Multiplication is associative since the composition of maps is.
Only commutativity and the existence of inverses need extra care.
For two slopes $\alpha,\beta$ we have the estimates
$$
n\alpha(\beta(n))=\alpha(n\beta(n))+E_1=\alpha(\beta(n)n)+E_1=\beta(n)\alpha(n)+E_2+E_1
$$
with $|E_1| \leq |n|S_{\alpha}$ and $|E_2|\leq |\beta(n)|S_{\alpha}\leq|n|(|\beta(1)+S_{\beta})S_{\alpha}$. It follows
$$
|\alpha\circ\beta(n)-\beta\circ\alpha(n)|
\leq S_{\alpha}(1+|\beta(1)|+S_{\beta})+S_{\beta}(1+|\alpha(1)|+S_{\alpha})
$$
showing that the slopes $\alpha \circ \beta$ and $\beta \circ \alpha$ are equivalent, hence the multiplication is
commutative.

Let $1$ be the real number represented by the identity map
$\Id_{\Z}:\Z \to \Z$. Clearly we have for a real number $x$ the properties
$1x=x1=x$, which makes $1$ into the unit element of the 
multiplication $\cdot$ 
of $\R$.

We now construct a right inverse for $x\in \R, x \not=0$, i.e. 
an element $y \in \R$ satisfying $xy=1$. Let $\alpha$ be a well adjusted 
representing slope for $x$. 
It follows, that for each $v \in \Z$ we may choose  $n_v \in \Z$ with
$|v-\alpha(n_v)| \leq |\alpha(1)|+1$. We define a map $\beta:\Z \to \Z$ by 
$\beta(v):=n_v$. 

We claim that the map $\beta$ is a slope. Indeed, for 
$v, w \in \Z$ we have 
$$
|\alpha(\beta(v+w)-\beta(v)-\beta(w))| = 
|\alpha(n_{v+w}-n_v-n_w)|\leq
$$
$$
|(v+w)-v-w|+2+3(|\alpha(1)|+1)= 3|\alpha(1)|+5.
$$
Since $\alpha$ takes each value only finitely many times, we conclude
that the set $\{\beta(v+w)-\beta(v)-\beta(w) \mid v,w \in \Z\}$ is finite.

For $v \in \Z$ we have 
$\alpha \circ \beta(v)=\alpha(n_v)$, so the 
slopes $\alpha \circ \beta$ and $\Id_{\Z}$ are 
equivalent, since $|v-\alpha(n_v)| \leq |\alpha(1)|+1$ holds. 
It follows that
$xy=1$.

The pair $(\R,<)$ is an order relation. First we prove that 
the relation $<$
is total. Let $x,y$ be real numbers 
represented by the slopes 
$\alpha$ and $\beta$. We consider the slope $\delta:=\alpha-\beta$, which
represents the number $x-y$. Let $\delta'$ be the well adjusted slope
equivalent to $\delta$ given by the concentration lemma. If $\delta'(n) \in
\{-1,0,1\}$ for all $n \in \Z$ then we have $x=y$. If $x\not=y$ we have 
for some $n \in \N$ either 
$\delta'(n)>1$, or $\delta'(n)<-1$. In the first case we have 
$x>y$ and in the second case we have $x<y$. The case $x=y$ excludes $x<y$ and
$x>y$. The cases $x<y$ and $x>y$ exclude each other. It remains 
only to prove
transitivity. Let $x,y,z$ be real numbers with 
$x>y$ and $y>z$, which are represented by the slopes                  
$\alpha, \beta$ and $\gamma$. Let 
$\delta_1, \delta_2$ be well adjusted slopes
equivalent to the slopes $\alpha-\beta$ and 
$\beta-\gamma$. So for some
$n \in \N_+$ and $m \in \N_+$ we have $\delta_1(n)>1$ and $\delta_2(m)>1$. It 
follows that $(\delta_1+\delta_2)(nm)>2$. The well adjusted slope
$\delta_{12}$ equivalent to $\delta_1+\delta_2$ will satisfy
$\delta_{12}(nm)>1$ and hence we have $x>z$. 

The quadruple $(\R,+,\cdot,<)$ is an ordered field. Let 
$x,y$ be reals satisfying 
$x<y$ and let $t$ be  real. We represent $x,y,t$ by the well adjusted slopes
$\alpha, \beta$ and $\tau$. Since $x<y$ there exists $b \in \N$ with 
$\alpha(bn) < \beta(bn)-n, n \in \N_+$. Hence, 
$\alpha(bn)+\tau(bn)<\beta(bn)+\tau(bn)-n, n \in \N_+$, 
showing the monotonicity property  for translations $x+t<y+t$. If $t>0$,
for some $d \in \N$ we have $\tau(dn)>n, n \in \N_+$, hence
$\tau(\alpha(bdn))< \tau(\beta(dbn))-n, n \in \N_+,$ showing the 
monotonicity property for the stretching $tx<ty$.	

We now prove the archimedean property. Let $a \in \R$ with $a>0$
and $A\in \R$  be given. We construct $N \in  \N$ such that
$Na > A$, as follows. We represent $a$  and $A$ by well adjusted slopes
$\lambda$ and $\Lambda$. Since $a > 0$ holds, we may choose $n \in \N_+$ 
with $\lambda(n)>1$. Then $\lambda(2n)>2$. Define 
$N:=1+\max\{\Lambda(2n),0\}$. Let $\kappa$ be the well adjusted slope
equivalent to the slope $N\lambda$. We have $\kappa(2n)> N\lambda(2n)-N >
2+\Lambda(2n)$. Hence $Na>A$.

Finally, we will establish the completeness. Let $D$ be a non-empty
subset in $\R$ bounded from above by $m \in \R$. So for $x \in D$ we have
the inequality $x \leq m$.
Let $\Delta$ be a set of well adjusted slopes representing the real numbers
in the set $D$. Let $\mu$ be a well adjusted slope representing $m$.
For every $n \in \N$ and every $\delta \in \Delta$ we have 
$\delta(n) < \mu(n)+2$. It follows that for $n \in \N_+$
the non-empty set $\{\delta(n) \mid \delta \in \Delta \}$ is 
bounded from above by $\mu(n)+2$. 
Let $\sigma : \Z \to \Z$ be the odd map defined by
$$
\sigma(n):=\max \{\delta(n) \mid \delta \in \Delta \}.
$$
We claim that the map $\sigma$ is a slope. Indeed, for 
$u \in \N_+$ let $\delta_u \in \Delta$ be a slope, which attains
at $u$ the value $\max\{\delta(u) \mid  \delta \in \Delta \}$. 
So, we have
$\delta_u(u)=\sigma(u).$ For $p, N \in \N_+$ put $q:=pN$. We compare 
$\delta_p, \delta_q$ at $p$ and $q$ as follows.
We have 
$$
\delta_q(q):N \leq \delta_p(p)+1
$$ 
since 
$|\delta_q(q):N-\delta_q(p)| \leq 1$ and $\delta_q(p) \leq \delta_p(p)$. 
We also have
$$
N\delta_p(p) \leq \delta_p(q)+N \leq \delta_q(q)+N.
$$ 
We conclude that for all $p ,N \in \N_+$
$$
|\delta_p(p)-\delta_{pN}(pN):N| \leq 1.
$$ 
Hence  for $n,m \in \N_+$, where we put 
$c:=\delta_{nm(n+m)}(nm(n+m))$, the following inequalities hold:
$$
|\sigma(n)-c:m(n+m)| \leq 1,
$$
$$
|\sigma(m)-c:n(n+m)| \leq 1,
$$
$$
|\sigma(n+m)-c:nm| \leq 1.
$$
For instance, the first inequality is obtained with $p=n, N=m(n+m), q=Np$ and
by comparing $\delta_p$ and $\delta_q$ at the point $p$. From 
$$
|c:nm-c:m(n+m)-c:n(n+m)| \leq 1
$$
it follows that for all $n,m \in \N_+$ we have
$$
|\sigma(n+m)-\sigma(n)-\sigma(m)| \leq 1+3=4,
$$
which proves our claim.

Let $s$ be the real number represented by the slope $\sigma$. 
For all $x \in D$ we have the inequality $x \leq s$, since for a slope
$\delta \in \Delta$ representing $x$ the inequalities
$$
\delta(n) \leq \delta_n(n)=\sigma(n), n \in \N_+,
$$
hold. So $s \in \R$ is an upper bound for $D$.

In order to prove that $s$ is the least upper bound of $D$, we show that 
no $t \in \R$ with $t<s$ is  an upper bound of $D$. Indeed, let
$\tau$ be a well adjusted slope for $t \in \R$ with $t<s$. There exists
$n \in \N_+$ with $\tau(n)<\sigma(n)-2$. Let $x$ in $D$ be represented by
$\delta_n$. We have $\delta_n(n) > \tau(n)+2$, hence $x>t$ and $t$ is not an
upper bound for $D$.

{\bf Remarks:} The above construction of the field of 
real numbers $(\R,+,\cdot)$ 
has as its starting point the additive group $(\Z,+)$. 
It is well known that the order relation $<$ on $\R$ is encoded in the
field structure of $\R$, namely for $x,y\in \R$ we have $x<y$ if and only if
$y-x=t^2$ holds for some $t \in \R \setminus \{0\}$. So we see that 
the real ordered field $(\R,+,\cdot,>)$
is constructed directly out of the additive group of integers. 

The group $\R/\Z$ appears as second bounded cohomology group $H^2_b(\Z,\Z)$ 
with coefficients $\Z$ of the group $\Z$. Bounded 
cohomology is defined by Michael Gromov in the seminal paper [G].
An element $\theta \in H^2_b(\Z,\Z)$ is by 
definition the class modulo the boundary $df$ of a bounded $1$-cochain
$f:\Z \to \Z$ of a bounded $2$-cocycle on $\Z$ with values in $\Z$. A bounded $2$-cocycle on 
$\Z$ with values in $\Z$
is a bounded map $\theta:\Z \times \Z \to \Z$ that 
satisfies $d\theta=0$. For example a real number $a$ 
defines a bounded $2$-cocycle $\theta_a(n,m):=[(n+m)a]-[na]-[ma]$. The map
$a \in \R \mapsto \theta_a$ induces an isomorphism from 
$\R/\Z$ to $H^2_b(\Z,\Z)$. We recall the formulae for differentials
in group cohomology: 
$d\theta(u,v,w):=\theta(u,v)-\theta(u+v,w)+\theta(u,v+w)-\theta(v,w)$ and 
$df(n,m):=f(n)-f(n+m)+f(m)$.
Our construction of the reals defines the additive group 
$(\R,+)$ as the quotient
$\R:=C^1_{db}(\Z,\Z)/C^1_b(\Z,\Z),$
with 
$$
C^1_{db}(\Z,\Z):=\{1{\it -cochains \,on\,}\, \Z 
{\it \,with\, values\, in \,}\Z {\it \,
and \,with\, bounded \,differential\,}\},
$$ 
$$
C^1_b(\Z,\Z):=  \{ {\it bounded}\, 1{\it -cochains\, on\,} 
\Z {\it \,with\, values\, in\,} \Z \}.
$$ 
The composition of maps in $C^1_{db}(\Z,\Z)=\{{\it slopes}\}$ induces the multiplication in $\R$. 

The encoding of the order relation on $\R$ in the field 
structure has far reaching consequences. For instance,
the fact that a line incidence preserving bijection
of the real projective plane is a projective transformation, is
such a consequence. 

The first rigorous definitions for real numbers were 
published independently in  $1872$ 
by G. Cantor, E. Heine, Ch. M\'eray, and R. Dedekind. The 
rigorous definition of convergence of
a sequence of numbers was given by d'Alembert in $1765$  and by Cauchy in 
$1821$ without
having at the time a rigorous definition of real numbers. 
An exposition of the construction
of real numbers is given 
in the book ``Grundlagen der Analysis'' of Edmund Landau [La].
We recommend  reading the book ``Analysis by Its History'' by E. Hairer 
and G. Wanner [H-W].

{\center{\bf References}}

\noindent
[A]
Niels Henrik Abel, {\it Beweis der Unm\"oglichkeit, algebraische 
Gleichungen von h\"oheren Graden als dem vierten allgemein aufzul\"osen}, 
J. reine angew. Math. {\bf 1} (1826), 65-84.

\noindent
[D]
Heinrich D\"orrie,
{\it Triumph der Mathematik : hundert ber\"uhmte Probleme aus zwei 
Jahrtausenden mathematischer Kultur}, Hirt, Breslau 1933, {\sl Eng. Transl.:} 
David Antin (1958),
{\it 100 Great Problems of Elementary Mathematics: Their History and Solution},
{\sl reprinted:} Dover, New York 1965.

\noindent
[F]
Kurt von Fritz, {\it The Discovery of Incommensurability by 
Hippasus of Metapontum}, Annals of Math. {\bf 46} 2, (1945), 242--264.

\noindent
[H-V]
D. Hilbert and S. Cohn-Vossen, {\it Anschauliche Geometrie}, 
Grundlehren der mathematischen Wissenschaften, 
{\bf 37} Springer, Berlin 1932.

\noindent
[H-W] 
E. Hairer and G. Wanner, {\it Analysis by its History}, Undergraduate
Texts in Mathematics, Springer Verlag, New-York, 
Berlin and Heidelberg, 1995.

\noindent
[G] 
M. Gromov, 
{\it Volume and bounded cohomology}, 
Publ. Math. I.H.E.S., 
{\bf 56} (1982), 5--100, (1983).

\noindent
[La] 
E. Landau, {\it Grundlagen der Analysis}, Akad. Verlags\-gesell\-schaft, 
Leipzig, 1930.

\noindent
[Le]
Henri Lebesgue, 
{\it Lecons sur les constructions g\'eom\'etriques}, avec
Pr\'eface de Paul Montel, Gauthier-Villars,
Paris 1950, and \'Editions Jacques Gabay, Paris 1987.

\noindent
[M]
Eli Maor,
{\bf e} {\it The Story of a Number},
Princeton University Press, Princeton 1994.

\noindent
[P] 
H. Poincar\'e, 
{\it Sur les courbes d\'efinies par les \'equations diff\'erentielles},
J. de Math. Pures et Appl., quatri\`eme s\'erie,  {\bf 1} (1885), 167--244 [Oeuvres, tome I, 90--161].

\noindent
[S]
C. Skau, {\it The impossibility of solving the general nth degree 
equation algebraically when $n \geq 5$ : Abel's 
and Ruffini's proofs revisited}, Normat {\bf 38} 2, (1990), 53-84.

\noindent
[R]
Paolo Ruffini, {\it Teoria generale delle equazioni, in cui 
si dimostra impossibile 
la soluzione algebraica delle equazioni generali 
di grado superiore al quarto}, Bologna, 
Stamperia di S. Tommaso d'Aquino, 1799. - 2 v. ; 22 cm. 
Volume I: VIII, 206, [4] p. 
Volume 2: [2], 207-509, [7] p., 2 tav. f. t. 

\noindent
[V] 
H. Vogt, {\it Die Entdeckungsgeschichte des Irrationalen nach 
Plato und anderen Quellen des 4. Jahrhunderts}, 
Bibliotheca Mathematica {\bf 3} (1910), 97--155.

\noindent
University of Basle, Rheinsprung 21, CH-4051  Basel.
\end{document}